\documentclass[12pt]{article}
\usepackage{amssymb,latexsym}
\newtheorem{theo}{Theorem}[section]
\newtheorem{defi}[theo]{Definition}
\newtheorem{lemm}[theo]{Lemma}

\newtheorem{prop}[theo]{Proposition}
\newtheorem{conj}[theo]{Conjecture}

\def\gl{\mathfrak{gl}}
\def\g{\mathfrak{g}}
\def\h{\mathfrak{h}}
\def\n{\mathfrak{n}}
\def\ch{\mathop{\rm ch}\nolimits}
\def\al{\alpha}
\def\be{\beta}
\def\ga{\gamma}
\def\de{\delta}
\def\ep{\epsilon}  
\def\varep{\varepsilon}

\def\th{\theta}
\def\la{\lambda}
\def\rh{\rho}
\def\si{\sigma}

\def\De{\Delta}

\def\beq{\begin{equation}}
\def\eeq{\end{equation}}
\def\bea{\begin{eqnarray}}
\def\eea{\end{eqnarray}}
\def\beas{\begin{eqnarray*}}
\def\eeas{\end{eqnarray*}}

\def\mybox{\hfill\llap{$\Box$}}
\begin{document}
\begin{center}
{\large \bf Characters and composition factor multiplicities}\\[3mm]
{\large \bf for the Lie superalgebra $\gl(m/n)$}\\[2cm]
{\sc J.\ Van der Jeugt}\footnote{Research Associate of the Fund
for Scientific Research -- Flanders (Belgium)} \\
Department of Applied Mathematics and Computer Science,\\
University of Ghent, Krijgslaan 281-S9, B-9000 Gent, Belgium\\[3mm]
{\sc R.B.\ Zhang}\\
Department of Pure Mathematics,\\
University of Adelaide, Adelaide, Australia
\end{center}

\begin{abstract}
The multiplicities $a_{\la,\mu}$ of simple modules
$L_\mu$ in the composition series of Kac modules $V_\la$ for the
Lie superalgebra $\gl(m/n)$ were described by Serganova, leading to her 
solution of the character problem for $\gl(m/n)$. In Serganova's algorithm
all $\mu$ with nonzero $a_{\la,\mu}$ are determined for a given $\la$;
this algorithm turns out to be rather complicated. 
In this Letter a simple rule is conjectured to find all nonzero 
$a_{\la,\mu}$ for any given weight $\mu$. In particular, we claim that
for an $r$-fold atypical weight $\mu$ there are $2^r$ distinct weights
$\la$ such that $a_{\la,\mu}=1$, and $a_{\la,\mu}=0$ for all other 
weights $\la$.
Some related properties on the
multiplicities $a_{\la,\mu}$ are proved, and arguments in favour of
our main conjecture are given. Finally, an extension of the conjecture
describing the inverse of the matrix of Kazhdan-Lusztig polynomials 
is discussed.
\end{abstract}

\vskip 1cm
{\small
\noindent 1991 Mathematics Subject Classification~: 17A70, 17B70.\\[1mm]
\noindent Keywords~: character formula, Lie superalgebra $\gl(m/n)$, Kac module,
composition factors, Kazhdan-Lusztig polynomial.
}
\vskip 1cm
\newpage
%
\section{Introduction}  \label{sec-intro}
\setcounter{equation}{0}

Shortly after the classification of finite-dimensional simple Lie 
superalgebras~\cite{Kac1,Scheunert}, Kac considered the problem of 
classifying all finite-dimensional simple modules (i.e.\ irreducible
representations) of the basic classical Lie superalgebras~\cite{Kac2}.
For a subclass of these simple modules, known as ``typical'' modules,
Kac derived a character formula closely analogous to the Weyl character
formula for simple modules of simple Lie algebras~\cite{Kac2}. The problem
of obtaining a character formula for the remaining ``atypical'' modules
has been the subject of intensive investigation, both in the mathematics
and physics literature. Several partial solutions to this problem were
given, e.g.\ for covariant or contravariant tensor 
representations~\cite{DJ,BR}, for so-called generic 
representations~\cite{PS}, for singly 
atypical representations~\cite{BL,VHKT2,VanderJeugt}, or for tame
representations~\cite{KW}. Only recently a solution for the characters
of simple $\gl(m/n)$ modules has been given by 
Serganova~\cite{Serganova2} (with partial results in~\cite{Serganova1}).
Thus the characters of all simple modules of type~I Lie superalgebras
are now in principle known, since for $C(n)$ the problem was already 
solved in~\cite{VanderJeugt}.
For the series of Lie superalgebras $\mathfrak{q}(n)$, the solution
of the character problem has been announced in~\cite{PS1}.

Let us briefly describe some aspects of Serganova's solution, and the
main results of the present paper.

The Lie superalgebra $\g=\gl(m/n)$ has a $\mathbb{Z}$-grading 
$\g=\g_{-1} \oplus \g_0 \oplus \g_{+1}$ consistent with the 
$\mathbb{Z}_2$-grading $\g=\g_{\bar 0} \oplus \g_{\bar 1}$, in
particular $\g_{\bar 0} = \g_0 = \gl(m) \oplus \gl(n)$. A 
finite-dimensional $\g_0$ module $L_\la(\g_0)$ with highest weight
$\la$ is turned into a $\g_0 \oplus \g_{+1}$ module by trivial $\g_{+1}$
action, and then the Kac module~\cite{Kac1,Kac2} is the induced module
\[
V_\la = {\cal U}(\g) \otimes_{ {\cal U}(\g_0 \oplus \g_{+1}) }
L_\la(\g_0\oplus \g_{+1}).
\]
This is a finite-dimensional module which is simple when $\la$ is 
typical, and indecomposable when $\la$ is atypical. The characters of
Kac modules are easy to determine. Any simple finite-dimensional
$\g$ module $L_\la = L_\la (\g)$ is the quotient of the Kac module
$V_\la$ by its unique maximal submodule. Since $V_\la$ is
indecomposable for atypical $\la$ one cannot decompose it into
irreducibles, but on the other hand $V_\la$ has a Jordan-H\"older 
composition series. Let $a_{\la,\mu}=[V_\la : L_\mu]$ be the multiplicity
of $L_\mu$ in the composition series of $V_\la$. Then
\[
\ch V_\la = \sum_\mu a_{\la,\mu} \ch L_\mu.
\]
With the canonical partial order on the set of weights, the matrix
$A=(a_{\la,\mu})$ is lower triangular (the diagonal consisting of 1's).
On the other hand, the character of $L_\la$ can be written as an
infinite ``alternating sum'' of  characters of Kac 
modules~\cite{Serganova2}, thus
\[
\ch L_\la = \sum_\mu b_{\la,\mu} \ch V_\mu.
\]
The matrix $B=(b_{\la,\mu})$ is also lower triangular, and $A$ and $B$
are inverses to each other. In a recent paper of
Serganova~\cite{Serganova2}, it is shown that
the coefficients $b_{\la,\mu}$ are equal to the value of Kazhdan-Lusztig
polynomials $K_{\la,\mu}(q)$ for $q=-1$, and an algorithm is given
for evaluating $b_{\la,\mu}=K_{\la,\mu}(-1)$ using induction on dimension
and the embedding $\gl(1)\oplus \gl(m-1/n) \subset \gl(m/n)$.
Serganova's work offers a principal solution to an outstanding problem 
which has been open for almost 20 years.
However, there remain some important open questions related to
the character problem for $\gl(m/n)$. 
First, the methods of~\cite{Serganova2} offer a direct solution to finding
the multiplicities $a_{\la,\mu}$, but to find the actual character
$\ch L_\la$ (i.e.\ the coefficients $b_{\la,\mu}$) the method is rather
indirect in the sense that it depends on an formal inversion process (see
also~\cite[Remarks~2.4, 2.5 and Example~2.6]{Serganova2}).
For example, it would be difficult to extract 
from~\cite[Theorem~2.3]{Serganova2} what the lowest weight $\mu$ is
in the decomposition of $L_\la$ with respect to the even subalgebra, i.e.
\[
L_\la (\g) = \bigoplus_\mu c_{\la,\mu} L_\mu (\g_0).
\]
Secondly, the algorithm presented in~\cite[Section~2]{Serganova2},
though straighforward, is not easy to apply, and does not shed much
light on the properties of the coefficients $a_{\la,\mu}$. Basically,
the algorithm starts with a given $\la$, and allows one in various steps
to determine the weights $\mu$ for which $a_{\la,\mu}$ is possibly not
zero; in general, however, many cancellations take place at the end
of the calculation and after this one finishes with the actual $\mu$'s
with non-zero $a_{\la,\mu}$.

The goal of the present Letter is to announce some results on
the multiplicities $a_{\la,\mu}$. In particular, we claim that
$a_{\la,\mu}$ is either 0 or 1. Moreover, we explain how the structure
of the matrix $A=(a_{\la,\mu})$, which seems extremely complicated 
when rows are considered, is in fact very simple when concentrating
on columns. In particular, we give an easy algorithm for finding 
the non-zero $a_{\la,\mu}$. Our algorithm is opposite in the sense that
for given $\mu$ all $\la$ are determined for which $a_{\la,\mu}\ne 0$,
and in such case $a_{\la,\mu}=1$; in all other cases $a_{\la,\mu}=0$.
We should immediately add that at the moment we have a proof of 
this property only in some cases, but our proof is not valid for
all cases. However, with all the data deduced by means of Serganova's
algorithm we feel safe that it is always valid. 
Based on this observation, we also present a conjecture giving the
inverse matrix $A_q=(a_{\la,\mu}(q))$ of the matrix of Kazhdan-Lusztig
polynomials $K_q=(K_{\la,\mu}(q))$. A strong argument in favour of 
this conjecture is presented, and some properties of Kazhdan-Lusztig
polynomials are derived.

The emphasis of this paper is on announcing these results, and on
proving some related properties. Some of the proofs involve
combinatorial arguments; here we have chosen to be as concise as
possible, so that the introduction of many combinatorial 
quantities~\cite{HKV} related to $\gl(m/n)$ can be avoided.

\section{Notation and definitions}  \label{sec-notation}
\setcounter{equation}{0}

Let $\g=\gl(m/n)$, $\h\subset\g$ its Cartan subalgebra, and 
$\g=\g_{-1} \oplus \g_0 \oplus \g_{+1}$ the consistent  
$\mathbb{Z}$-grading. Note that $\g_{0}=\gl(m)\oplus \gl(n)$, and
put $\g^+=\g_0\oplus \g_{+1}$ and $\g^-=\g_0\oplus \g_{-1}$. The
dual space $\h^*$ has a natural basis $\{\ep_1,\ldots,\ep_m, \de_1,
\ldots, \de_n\}$, and the roots of $\g$ can be expressed in terms
of this basis. Let $\De$ be the set of all roots, $\De_0$ the set of
even roots, and $\De_1$ the set of odd roots. One can choose a set of
simple roots (or, equivalently, a triangular decomposition), but note
that contrary to the case of simple Lie algebras not all such choices
are equivalent. The distinguished choice for a triangular decomposition
$\g=\n^-\oplus\h\oplus\n^+$ is such that $\g_{+1}\subset \n^+$ and
$\g_{-1}\subset \n^-$. Then $\h\oplus\n^+$ is the corresponding
distinguished Borel subalgebra, and ${\De_+}$ the set
of positive roots. For this choice we have explicitly~:
\beas
\De_{0,+} &=& \{\ep_i-\ep_j | 1\leq i < j \leq m \} \cup
\{\de_i-\de_j | 1\leq i < j \leq n \} , \\
\De_{1,+} &=& \{ \be_{ij}=\ep_i-\de_j | 1\leq i \leq m,\ 1\leq j \leq n\},
\eeas
and the corresponding set of simple roots is given by 
\[
\{ \ep_1-\ep_2,\ldots,\ep_{m-1}-\ep_m, \ep_m-\de_1, \de_1-\de_2,
\ldots,\de_{n-1}-\de_n\}.
\]
Thus in the distinguished basis there is only one simple root which is odd.
As usual, we put
\[
\rho_0 = {1\over 2} \bigl( \sum_{\al\in\De_{0,+}} \al\bigr),\qquad
\rho_1 = {1\over 2} \bigl( \sum_{\al\in\De_{1,+}} \al \bigr),\qquad
\rh=\rh_0-\rh_1.
\]
There is a symmetric form $(\ ,\ )$ on $\h^*$ induced by the invariant
symmetric form on $\g$, and in the natural basis it takes the form
$(\ep_i,\ep_j)=\de_{ij}$, $(\ep_i,\de_j)=0$ and 
$(\de_i,\de_j)=-\de_{ij}$. 
A weight $\la\in\h^*$ is integral if $(\la,\al)\in\mathbb{Z}$ for all
even roots $\al$, and it is dominant if 
$2(\la,\al)/(\al,\al)$ is a nonnegative
integer for all positive even roots $\al$. Following 
Sergavona~\cite{Serganova2} we say that $\la$ is regular if it is not
on a wall in any Weyl chamber, i.e.\ if $(\la,\al)\ne 0$ for every
$\al$ in $\De_0$. The set of integral weights is denoted by $P$,
the set of integral dominant weights by $P^+$. 
With our choice of positive roots, the weights $P$ are partially
ordered by $\la\leq\mu$ iff $\mu-\la=\sum {k_\al}\al$ where $\al\in{\De_+}$
and $k_\al$ nonnegative integers. In the standard $\ep$-$\de$-basis,
an integral weight $\la\in P$ is written as
\beas
\la&=& \la_1 \ep_1 + \cdots + \la_m\ep_m + \la_1'\de_1 + \cdots + 
\la_n'\de_n,\\
&=& (\la_1,\cdots,\la_m; \la_1',\cdots,\la_n'),
\eeas
where $\la_i-\la_{i+1}\in\mathbb{Z}$ and $\la_i'-\la_{i+1}'\in\mathbb{Z}$;
$\la\in P^+$ if moreover $\la_1\geq\la_2\geq\cdots \la_m$ and
$\la_1'\geq\la_2'\geq\cdots \la_n'$. 

The Weyl group of $\g$ is the Weyl group $W$ of $\g_0$, hence it is
the direct product of symmetric groups $S_m\times S_n$. For $w\in W$,
we denote by $\varep(w)$ its signature. The dot action is defined as
$w\cdot\la = w(\la+\rh)-\rh$, for $w\in W$ and $\la\in\h^*$. For $\la
\in P$, $W\la$ has a unique element in $P^+$ (i.e.\ in the dominant
Weyl chamber), and this is denoted by $d(\la)$. If $d(\la+\rh)-\rh \in
P^+$, then we define ${\dot d}(\la)=d(\la+\rh)-\rh$; if 
$d(\la+\rh)-\rh \not\in P^+$, then ${\dot d}(\la)$ is said to be 
undefined.

For $\la$ integral dominant, let $L_\la(\g_0)$ denote the 
finite-dimensional irreducible $\g_0$ module with highest weight $\la$.
This can be extended to a $\g^+$ module $L_\la(\g^+)$
by letting $\g_{+1}$ act trivially
on the elements of $L_\la(\g_0)$. The Kac module $V_\la$ is then the
finite-dimensional module ${\cal U}(\g) \otimes_{ {\cal U}(\g^+) }
L_\la(\g^+)$. This is in general indecomposable, and the irreducible
$\g$ module obtained by quotienting $V_\la$ by its maximal submodule
is denoted by $L_\la$ or $L_\la(\g)$. All these modules $V$ are 
$\h$-diagonalizable with weight decomposition $V=\oplus_\mu V(\mu)$, 
and the character is defined to be $\ch V = \sum_\mu \dim V(\mu) \, 
e^\mu$, where $e^\mu$ is the formal exponential.

Let $\la\in P^+$, then $\la$ (resp.\ $V_\la$ and $L_\la$) is said
to be typical if $(\la+\rh,\al)\ne 0$ for every $\al\in\De_{1,+}$. 
Otherwise $\la$ (resp.\ $V_\la$ and $L_\la$) is said to be atypical.
The number $r=\# \la$ of elements $\al\in\De_{1,+}$ for which 
$(\la+\rh,\al)= 0$ is called the degree of atypicality. This definition
depends upon the choice of $\De_+$, but one can easily extend it such
that the degree of atypicality of a simple $\g$ module $V$ does not
depend upon this choice (see~\cite[Corollary~3.1]{KW}).
If $\la$ is typical, $L_\la = V_\la$, and the character is easy to
write down~:
\[
\ch V_\la = \prod_{\al\in\De_{1,+}} (1+e^{-\al}) \,\ch L_\la(\g_0).
\]
If $\la$ is atypical then $V_\la$ is indecomposable, and we denote
by $a_{\la,\mu}=[V_\la : L_\mu]$ the multiplicity
of $L_\mu$ in the composition series of $V_\la$, thus
$\ch V_\la = \sum_\mu a_{\la,\mu} \ch L_\mu$.
The multiplicities $a_{\la,\mu}$ can be nonzero only when $\mu \leq \la$
with respect to the partial order. Moreover $a_{\la,\la}=1$ and
if $\#\mu\ne\#\la$ then $a_{\la,\mu}= 0$ , see~\cite{Serganova2}.

\section{Preliminaries}  \label{sec-prelim}
\setcounter{equation}{0}

In this section (and the following one) $\mu$ is an integral
dominant $r$-fold atypical weight ($\#\mu=r$), and the $r$ elements $\al$
of $\De_{1,+}$ for which $(\mu+\rh,\al)=0$ are denoted by
$\{\ga_1,\ldots,\ga_r\}$. Moreover, they are ordered in such a way that
\[
\ga_1< \ga_2 < \cdots < \ga_r.
\]

\begin{lemm}
Let $\mu$ be integral dominant and $\al\in\De_{1,+}$ such that
$(\mu+\rh,\al)=0$. Then there exists a unique maximal subset
$\{ \al=\al_0, \al_1, \ldots ,\al_k\}$ of $\De_{1,+}$ with every
$\al_i > \al$ (for $i>0$) such that
\[
 (\mu+\rh,\al_0)=0, (\mu+\al_0+\rh,\al_1)=0, \ldots,
(\mu+\al_0+\cdots +\al_{k-1}+\rh,\al_k)=0.
\]
Moreover, $\la=\mu+\al_0+\cdots +\al_{k-1}+\al_k$ is integral dominant.
\label{lemm1}
\end{lemm}
\noindent
The set is maximal in the sense that any element $\al_{k+1}>\al$ for
which $(\la+\rh,\al_{k+1})=0$ belongs already to the set.

\noindent
{\em Proof.} The proof is combinatorial, and we do not give it here.
It follows exactly the same arguments as in~\cite{PS0}, or as 
in~\cite[Definition~6.1 and Theorem~6.2]{HKV}. \mybox

\begin{defi}
Let $\mu$ be integral dominant and $r$-fold atypical 
with respect to the roots $\ga_1< \ga_2 < \cdots < \ga_r$ ($\ga_i \in
\De_{1,+}$). For each $\ga_i$, let $\De(\ga_i)$ denote the subset
of $\De_{1,+}$ defined by Lemma~\ref{lemm1}.
\label{defi1}
\end{defi}

It is easy to deduce that
\[
\De(\ga_1) \supset \De(\ga_2) \supset \cdots \supset \De(\ga_r).
\]

\addtocounter{theo}{1}
\noindent {\bf Example \arabic{section}.\arabic{theo}}
Let $\g=\gl(4/5)$, and $\mu=(2,1,0,0;0,-2,-2,-2,-2)$. The numbers
$(\mu+\rh,\be_{i,j})$ ($\be_{i,j}=\ep_i-\de_j$) are given here in
the atypicality matrix~\cite{VHKT1,VHKT2} $A(\mu)$~:
\[
A(\mu)=\left( \begin{array}{rrrrr}
5&2&1&0&-1\\ 3&0&-1&-2&-3 \\ 1&-2&-3&-4&-5\\ 0&-3&-4&-5&-6 \end{array}
\right).
\]
Thus $\#\mu=3$, $\ga_1=\be_{4,1}$, $\ga_2=\be_{2,2}$, and 
$\ga_3=\be_{1,4}$. One can verify that~:
\beas
\De(\ga_3)&=& \{ \be_{1,4}, \be_{1,5} \},\\
\De(\ga_2)&=& \{ \be_{2,2}, \be_{2,3}, \be_{2,4}, \be_{2,5}, 
\be_{1,3}, \be_{1,4}, \be_{1,5} \},\\
\De(\ga_1)&=& \{ \be_{4,1}, \be_{3,1}, 
\be_{2,2}, \be_{2,3}, \be_{2,4}, \be_{2,5}, 
\be_{1,3}, \be_{1,4}, \be_{1,5} \}.
\eeas

Since $\De(\ga_i)\supset \De(\ga_{i+1})$ we can consider their 
differences.
\begin{defi}
Let $\mu$ be integral dominant and $r$-fold atypical 
with respect to the roots $\ga_1< \ga_2 < \cdots < \ga_r$, and
$\De(\ga_i)$ as in Definition~\ref{defi1}. Then
$\nabla(\ga_i)=\De(\ga_i) \backslash \De(\ga_{i+1})$ (with
$\nabla(\ga_r)=\De(\ga_r)$). Denote by $k_i=\#\nabla(\ga_i)$
the number of elements in $\nabla(\ga_i)$. 
Two such roots $\ga_i$ and $\ga_j$ are said to be disconnected for $\mu$
if $\nabla(\ga_i) \perp \nabla(\ga_j)$ (orthogonal with respect to
the symmetric form $(\ ,\ )$ on $\h^*$); otherwise they are connected
for $\mu$.
\label{definabla}
\end{defi}

In example~3.3, we have 
\beas
\nabla(\ga_3)&=& \{ \be_{1,4}, \be_{1,5} \},\\
\nabla(\ga_2)&=& \{ \be_{2,2}, \be_{2,3}, \be_{2,4}, \be_{2,5}, 
\be_{1,3} \},\\
\nabla(\ga_1)&=& \{ \be_{4,1}, \be_{3,1} \},
\eeas
hence $(k_1,k_2,k_3)=(2,5,2)$; $\ga_1$ and $\ga_2$ (or $\ga_3$) are
disconnected, whereas $\ga_2$ and $\ga_3$ are connected for $\mu$.

\begin{prop}
Let $\mu$ be integral dominant and $r$-fold atypical 
with respect to the roots $\ga_1< \ga_2 < \cdots < \ga_r$ ($\ga_i \in
\De_{1,+}$). Consider the decomposition of $L_\mu$ with respect to the
even subalgebra~:
\[
L_\mu=L_\mu (\g) = \bigoplus_\nu c_{\mu,\nu} L_\nu (\g_0).
\]
The set $S=\{ \nu\in P^+ | c_{\mu,\nu}\ne 0\}$ contains a unique smallest
element $\mu_0$ (with respect to the partial order), 
$c_{\mu,\mu_0}=1$, and
\[
\mu_0 = \mu - \sum_{\al \in \De_{1,+} \backslash \De(\ga_1)} \al.
\]
\label{prop-mu0}
\end{prop}

\noindent {\em Proof.} 
We need the notion of reflection with respect to an odd simple 
root~\cite{PS0,KW}. For a given subset of positive roots ${\De_+}$
of $\De$ and an odd simple root $\al$, one may construct a new subset
of positive roots $\De_+'$ by
\[
\De_+'=(\De_+ \cup \{ -\al\} ) \backslash \{\al\}.
\]
The set $\De_+'$ is said to be obtained from $\De_+$ by a simple 
$\al$-reflection. Note that the set of even positive roots $\De_{0,+}$
remains unchanged. In that case, the new $\rh$ is given by
\[
\rh' = \rh + \al.
\]
When $V$ is a finite-dimensional simple $\g$ module, and $\De_+$ is
a fixed set of positive roots, then $V$ has a highest weight $\la$ with
respect to $\De_+$. If $\De_+'$ is obtained from $\De_+$ by a simple
$\al$-reflection ($\al$ odd), then the highest weight $\la'$ of $V$
with respect to $\De_+'$ is given by~\cite[(3.1)]{KW}~:
\[
\begin{array}{lcl}
\la' = \la-\al & \hbox{if} & (\la | \al) \ne 0 \\
\la' = \la     & \hbox{if} & (\la | \al) = 0.
\end{array}
\]
Moreover, one has
\[
\begin{array}{lcl}
\la'+\rh' = \la+\rh & \hbox{if} & (\la+\rh | \al) \ne 0 \\
\la'+\rh' = \la+\rh+\al     & \hbox{if} & (\la+\rh | \al) = 0.
\end{array}
\]

\noindent
Now we apply this to $V=L_\mu(\g)$, starting with the distinguished
positive set of roots $\De_+$ given in Section~\ref{sec-notation},
thus the highest weight with respect to $\De_+$ is $\mu$.
The only odd simple root is $\al=\be_{m,1}=\ep_m-\de_1$; applying the
simple $\al$-reflection leads to a new set of positive roots $\De_+'$
and the new highest weight $\mu'$ is given by $\mu'=\mu-\al$ if
$(\mu+\rh|\al)\ne 0$ or by $\mu'=\mu$ if $(\mu+\rh| \al)=0$
(for the current $\al$ we have that $(\rh|\al)=0$).
Next, apply the simple $\al'$-reflection with respect to $\al'=\be_{m,2} =
\ep_m-\de_2$ (which is indeed an odd simple root in 
$\De_+'$) to obtain the new set $\De_+''$.
We have $(\mu'|\al')=(\mu'+\rh'|\al')$. So
\begin{itemize}
\item if $(\mu+\rh | \al) \ne 0$ then
 \begin{itemize}
 \item if $(\mu+\rh | \al') =(\mu'+\rh'|\al')=(\mu'|\al')\ne 0$ then
       $\mu''=\mu'-\al'=\mu-\al-\al'$,
 \item if $(\mu+\rh|\al')=0$ then $\mu''=\mu'=\mu-\al$;
 \end{itemize}
\item if $(\mu+\rh | \al) = 0$ then
 \begin{itemize}
 \item if $(\mu+\al+\rh | \al') =(\mu'+\rh'|\al')=(\mu'|\al')\ne 0$ then
       $\mu''=\mu'-\al'=\mu-\al'$,
 \item if $(\mu+\al+\rh|\al')=0$ then $\mu''=\mu'=\mu$.
 \end{itemize}
\end{itemize}
So far we have applied only two simple $\al$-reflections, namely first
with respect to $\be_{m,1}$ and then with respect to $\be_{m,2}$. 
Now continue applying such reflections, with respect to (in this order)
\[
\be_{m,3},\ldots,\be_{m,n},\be_{m-1,1},\be_{m-1,2},\ldots,\be_{m-1,n},
\ldots, \ldots, \be_{1,1},\be_{1,2}, \ldots, \be_{1,n}.
\]
Using lemma~\ref{lemm1} and definition~\ref{defi1},
it is an easy combinatorial exercise to see that with respect to the
final set of positive roots, the highest weight of our module $V$ is
given by $\mu_0=\mu - \sum_{\al \in \De_{1,+} \backslash \De(\ga_1)} \al$,
since in this process an odd root $\al$ will be subtracted only when it
is not in $\De(\ga_1)$. The final set of positive roots is given
by 
\[
\De_{f,+}=\De_{0,+} \cup \{ -\be_{i,j}=\de_j - \ep_i | 1\leq i\leq m,
 \ 1\leq j \leq n\},
\]
where $\De_{0,+}$ remains unchanged. 
But at every stage in the process, the new highest weight $\nu$ is
unique and satisfies $c_{\mu,\nu}=1$. The final weight $\mu_0$ is the
highest weight of $V$ with respect to $\De_{f,+}$, so it must be the 
smallest element in $S$. \mybox

\noindent
In example~3.3 we have $\mu_0=(0,0,-4,-4; 2,1,0,0,0)$. This is the
lowest $\g_0$ highest weight in the decomposition of $L_\mu$ with respect
to the even subalgebra $\g_0$.

\begin{prop}
Let $\mu$ be integral dominant and $r$-fold atypical 
with respect to the roots $\ga_1< \ga_2 < \cdots < \ga_r$ ($\ga_i \in
\De_{1,+}$), and
\[
\mu_0 = \mu - \sum_{\al \in \De_{1,+} \backslash \De(\ga_1)} \al.
\]
Let $k_i=\#\nabla(\ga_i)$ (see Definition~\ref{definabla}). Then
\[
{\dot d}(\mu+\sum_{i=1}^r k_i \ga_i) = \mu_0 + 2\rh_1.
\]
\end{prop}

{\em Proof.} The proof is combinatorial and uses a number of notions
defined in~\cite{HKV}. With the weight $\mu$ there corresponds a
composite Young diagram, specified by the Young diagrams of the 
two parts of $\mu$ associated to $\gl(m)$ and $\gl(n)$ respectively.
The addition of a coordinated boundary strip of length $k$ to $\mu$,
starting at the position of $\ga_i$, was 
discussed in~\cite{HKV}, and the corresponding weight is given
by ${\dot d}(\mu+k_i \ga_i)$ if the resulting composite diagram is 
standard (if the resulting composite diagram is not standard,
${\dot d}$ is undefined on $\mu+k_i \ga_i$). 
Using the arguments of~\cite[Theorem~6.2]{HKV} or 
of~\cite[Lemma~6.7]{VHKT2} one can deduce that
\[
{\dot d}(\mu+k_r\ga_r) = \mu+\sum_{\al\in\nabla(\ga_r)} \al \qquad\in P^+,
\]
and then
\beas
{\dot d}(\mu+\sum_{i=1}^r k_i\ga_i) &=& 
\mu+\sum_{\al\in\nabla(\ga_1)} \al +\sum_{\al\in\nabla(\ga_2)} \al + \cdots
+\sum_{\al\in\nabla(\ga_r)} \al \\
&=& \mu + \sum_{\al\in\De(\ga_1)} \al.
\eeas
Using $\mu_0$ and the definition of $\rh_1$, the final result 
follows. \mybox

\noindent
In example~3.3, 
\[
\mu+\sum_{i=1}^r k_i\ga_i=\mu+2\be_{4,1}+5\be_{2,2}+
2\be_{1,4} = (4,6,0,2;-2,-7,-2,-4,-2). 
\]
Applying ${\dot d}$ to this weight
gives $(5,5,1,1;-2,-3,-4,-4,-4)$, which is indeed $\mu_0+2\rh_1$.

\section{Main results}  \label{sec-result}
\setcounter{equation}{0}

Suppose again that $\mu$ is integral dominant and $r$-fold atypical 
with respect to the roots $\ga_1< \ga_2 < \cdots < \ga_r$ ($\ga_i \in
\De_{1,+}$). The Kac module $V_{\mu_0+2\rh_1}$, with highest weight 
$\mu_0+2\rh_1$, is completely reducible with respect to $\g_0$. From
the structure of Kac modules it follows that in this reduction 
$V_{\mu_0+2\rh_1}$ has a unique $\g_0$ component with highest weight 
$\mu_0$. In fact, the highest weight vector $v_{\mu_0}$ of this 
component is obtained by applying the product of all negative odd root
vectors to the highest weight vector $v_{\mu+2\rh_1}$ of the Kac 
module~\cite{Kac2,VHKT2}. This vector $v_{\mu_0}$ is contained in
every submodule of $V_{\mu_0+2\rh_1}$. Then it follows from 
proposition~{\ref{prop-mu0}} 
that every $\g$ submodule of $V_{\mu_0+2\rh_1}$
must also contain $L_\mu$, in other words $L_\mu$ is the smallest
$\g$ submodule in the Kac module $V_{\mu_0+2\rh_1}$.

Thus, starting from $\mu$ one determines the weight 
$\la={\dot d}(\mu+\sum_{i=1}^r k_i\ga_i)$. Then $V_\la$ has $L_\mu$
as a composition factor, $a_{\la,\mu}=1$, and $L_\mu$ is the
smallest $\g$ submodule in $V_\la$. Moreover, for every $\nu > \la$,
we have that $a_{\nu,\mu}=0$. Indeed, $V_\nu$ cannot have $L_\mu$
as a composition factor since the smallest $\g_0$ highest weight in
$V_\nu$ is $\nu-2\rh_1$, and $\nu-2\rh_1 > \mu_0$.

Hence, all weights $\la$ for which $a_{\la,\mu}\ne 0$ must lie
between $\mu \leq \la \leq {\dot d}(\mu+\sum_{i=1}^r k_i\ga_i)$.
The following is our main result~: we give a simple expression for
those $\la$ with $a_{\la,\mu}\ne 0$~:

\begin{conj}
Let $\mu$ be integral dominant and $r$-fold atypical 
with respect to the roots $\ga_1< \ga_2 < \cdots < \ga_r$ ($\ga_i \in
\De_{1,+}$), and $k_i=\#\nabla(\ga_i)$ (see definition~\ref{definabla}). 
For $\th=(\th_1,\ldots,\th_r)\in \{0,1\}^r$, consider
\[
\la_\th = {\dot d}(\mu_\th) = {\dot d}(\mu+\sum_{i=1}^r \th_i k_i \ga_i).
\]
Then $a_{\la_\th,\mu}=1$ for each of the $2^r$ integral dominant
weights $\la_\th$, and $a_{\la,\mu}=0$ elsewhere.
\label{conj1}
\end{conj}

There are a number of cases in which the conjecture can be proved.
For $r=0$ it is trivial; for $r=1$ it follows from the results 
of~\cite{VHKT2}.
For generic weights $\mu$ it can be deduced from~\cite{PS0,PS}. 
When all roots $\ga_i$ are disconnected for $\mu$, the conjecture can
be proved by induction on $r$. But for the general case we do not
have a proof so far. For every given $\la$, 
$a_{\la,\mu}$ can be calculated from Serganova's 
algorithm~\cite{Serganova2}, and thus we were able to successfully
verify this conjecture for numerous examples. In Serganova's
approach, induction from $\gl(1)\oplus \gl(m-1/n)$ to $\gl(m/n)$ is
used, and in this setting the question of finding those $\mu$
with $a_{\la,\mu}\ne 0$ for given $\la$ is natural; but the question
of finding those $\la$ with $a_{\la,\mu}\ne 0$ for given $\mu$ 
is rather unnatural and we think it cannot be solved directly using
the same type of induction.

Finally, note that the above conjecture is closely related 
to~\cite[Conjecture~7.2]{HKV}.

\addtocounter{theo}{1}
\noindent {\bf Remark \arabic{section}.\arabic{theo}} \ 
Since BGG duality holds for $\gl(m/n)$~\cite[Theorem~2.7]{Zou}, it 
follows that $a_{\la,\mu}$ also describes the multiplicity of the Kac
module $V_\la$ in a Kac composition series of the indecomposable 
projective module $I_\mu$~\cite[\S 2]{Zou}. Thus according to 
Conjecture~\ref{conj1} the multiplicities of Kac modules in 
indecomposable projective modules are easier to describe than the
multiplicities of simple modules in Kac modules.

\addtocounter{theo}{1}
\noindent {\bf Example \arabic{section}.\arabic{theo}} \ 
Take the data from example~3.3, i.e.\
$\g=\gl(4/5)$, and\\ $\mu=(2,1,0,0;0,-2,-2,-2,-2)$. We know already
that $r=3$, $(\ga_1,\ga_2,\ga_3)=(\be_{4,1},\be_{2,2},\be_{1,4})$,
and $(k_1,k_2,k_3)=(2,5,2)$. It is easy to calculate the $\mu_\th$
and $\la_\th$~:
\[
\begin{array}{cll}
\th~: & \mu_\th~: & \la_\th~: \\[2mm]
(0,0,0) &  (2,1,0,0; 0,-2,-2,-2,-2)  &  (2,1,0,0; 0,-2,-2,-2,-2) \\
(1,0,0) &  (2,1,0,2;-2,-2,-2,-2,-2)  &  (2,1,1,1;-2,-2,-2,-2,-2) \\
(0,1,0) &  (2,6,0,0; 0,-7,-2,-2,-2)  &  (5,3,0,0; 0,-3,-3,-3,-4) \\
(0,0,1) &  (4,1,0,0; 0,-2,-2,-4,-2)  &  (4,1,0,0; 0,-2,-2,-3,-3) \\
(1,1,0) &  (2,6,0,2;-2,-7,-2,-2,-2)  &  (5,3,1,1;-2,-3,-3,-3,-4) \\
(1,0,1) &  (4,1,0,2;-2,-2,-2,-4,-2)  &  (4,1,1,1;-2,-2,-2,-3,-3) \\
(0,1,1) &  (4,6,0,0; 0,-7,-2,-4,-2)  &  (5,5,0,0; 0,-3,-4,-4,-4) \\
(1,1,1) &  (4,6,0,2;-2,-7,-2,-4,-2)  &  (5,5,1,1;-2,-3,-4,-4,-4) 
\end{array}
\]
So all Kac modules $V_{\la_\th}$ have $L_\mu$ as a composition factor,
with multiplicity 1, and $a_{\la,\mu}=0$ for all other $\la$.
The first weight is $\la_{(0,0,0)}=\mu$, and obviously $a_{\mu,\mu}=1$;
the last weight is $\la_{(1,1,1)}=\mu_0+2\rh_1$, and here we have proved
earlier that indeed $a_{\mu_0+2\rh_1,\mu}=1$.

We shall now consider a number of consequences of Conjecture~\ref{conj1}.
First, the matrix $A=(a_{\la,\mu})$ is now easy to determine,
column by column. With
$\#\mu$ denoting the degree of atypicality, we have
\[
\sum_\la a_{\la,\mu} = 2^{\#\mu}
\]
for every $\mu\in P^+$. Recall that the inverse matrix $B=(b_{\la,\mu})$
of the lower triangular matrix $A$ consists of the coefficients in the
character formula for $L_\la$~:
\[
\ch L_\la = \sum_\mu b_{\la,\mu} \ch V_\mu.
\]
These coefficients are equal to specializations of Kazhdan-Lusztig
polynomials $K_{\la,\mu}(q)$, i.e.\ $b_{\la,\mu}=K_{\la,\mu}(-1)$.
Kazhdan-Lusztig polynomials for $\gl(m/n)$ were defined by 
Serganova~\cite{Serganova1,Serganova2}. Consider the $i$th homology
$H_i(\g_{-1}; L_\la)$. This space has the structure of a $\g_0$ module,
and denote the multiplicity $[H_i(\g_{-1}; L_\la) : L_\mu(\g_0)]$
by $K_{\la,\mu}^i$. Then the Kazhdan-Lusztig polynomials are defined as
\[
K_{\la,\mu}(q) = \sum_{i=0}^\infty K_{\la,\mu}^i q^i.
\]
Since $H_i(\g_{-1}; L_\la)$ is a $\g_0$ quotient module of 
$\hbox{Sym}^i(\g_{-1}) 
\otimes L_\la$, the $K_{\la,\mu}(q)$'s are polynomials (and not infinite
series) in $q$.

In general $H_i(\g_{-1}; L_\la)$ is difficult to determine, except when
$\la=0$ because in that case
\[
H_i(\g_{-1}; L_0) = \hbox{Sym}^i (\g_{-1}).
\]
One can explicitly construct the decomposition with respect to $\g_0$
for these modules (assume $m\leq n$)~:
\[
\hbox{Sym}^i (\g_{-1}) = \bigoplus_\si 
L_{(-\si_m,\ldots,-\si_2,-\si_1;
\si_1,\si_2,\ldots,\si_m,0,\ldots,0)}(\g_0),
\]
where the sum is over those partitions $\si$ with $|\si|=i$ and with
$m$ parts, i.e.\ over 
all integers $\si_1\geq \si_2 \geq \cdots \geq \si_m \geq 0$ with
$\si_1+ \si_2 + \cdots + \si_m =i$. Thus
\[
K_{0, (-\si_m,\ldots,-\si_2,-\si_1;
\si_1,\si_2,\ldots,\si_m,0,\ldots,0)}(q) = 
q^{\si_1+ \si_2 + \cdots + \si_m}
\]
if $\si$ is a partition with $m$ parts, and $K_{0,\mu}(q)=0$ for all
other weights $\mu$. Observe that this result is in agreement with
the expansion~(8.8) of~\cite{VHKT1} for $q=-1$.

Since we have a simple way to calculate the matrix $A=(a_{\la,\mu})$,
being the inverse of $B=(b_{\la,\mu})=(K_{\la,\mu}(-1))$,
it would be interesting to see if also the inverse of 
$K_q=(K_{\la,\mu}(q))$,
say $A_q = (a_{\la,\mu}(q))$, can easily be determined. 
Here, we have the following~:

\begin{conj}
Let $\mu$ be integral dominant and $r$-fold atypical 
with respect to the roots $\ga_1< \ga_2 < \cdots < \ga_r$ ($\ga_i \in
\De_{1,+}$). Consider the $2^r$ integral dominant weights $\la_\th$,
$\th=(\th_1,\ldots,\th_r)\in \{0,1\}^r$, determined in 
conjecture~\ref{conj1}. Let
\[
a_{\la_\th,\mu}(q)= (-q)^{|\th|},\qquad |\th|=\sum_i \th_i,
\]
and $a_{\la,\mu}(q)=0$ for all other $\la$.
Then the inverse of the triangular matrix $A_q = ( a_{\la,\mu}(q) )$ 
is the matrix of Kazhdan-Lusztig polynomials $K_q=(K_{\la,\mu}(q) )$.
\label{conj2}
\end{conj}

We have no general proof, but only a number of consistency checks.
For this purpose, it is useful to consider in $A_q$ and $K_q$ the
submatrices corresponding to $r$-fold atypical weights (there is no
overlap in these submatrices, i.e.\ if $\#\la \ne \#\mu$ then
$a_{\la,\mu}(q)=0=K_{\la,\mu}(q)$.) For the submatrices corresponding
to typical weights, the conjecture obviously holds, since both these
submatrices are the identity matrix. For the submatrices
corresponding to singly atypical weights ($r=1$), it follows 
from~\cite{VHKT2} and~\cite{Zou} that the conjecture is true.
As an additional verification, we have considered the determination
of $K_{0,\mu}(q)$ by explicitly inverting the matrix $(a_{\la,\mu}(q))$
for some $\gl(m/n)$ with small values of $m$ and $n$~:

\addtocounter{theo}{1}
\noindent {\bf Example \arabic{section}.\arabic{theo}}
Let $\g=\gl(2/2)$, and consider the integral dominant
weights $\mu$ with $\#\mu=2$.
These are of the form $(x,y;-y,-x)$, $x$ and $y$ integers with
$x\geq y$. The nonzero multiplicities are given by
\begin{itemize}
\item if $x=y$ then $a_{(x,y;-y,-x),(x,y;-y,-x)}=1$, 
$a_{(x,y;-y,-x),(x,y-1;-y+1,-x)}=-q$, and
$a_{(x,y;-y,-x),(x-2,y-2;-y+2,-x+2)}=q^2$.
\item if $x=y+1$ then
$a_{(x,y;-y,-x),(x,y;-y,-x)}=1$,
$a_{(x,y;-y,-x),(x,y-1;-y+1,-x)}=-q$,
$a_{(x,y;-y,-x),(x-1,y;-y,-x+1)}=-q$,
$a_{(x,y;-y,-x),(x-2,y-1;-y+1,-x+2)}=-q$, and\\
$a_{(x,y;-y,-x),(x-1,y-1;-y+1,-x+1)}=q^2$.
\item if $x\geq y+2$ then
$a_{(x,y;-y,-x),(x,y;-y,-x)}=1$,
$a_{(x,y;-y,-x),(x-1,y;-y,-x+1)}=-q$,
$a_{(x,y;-y,-x),(x,y-1;-y+1,-x)}=-q$, and
$a_{(x,y;-y,-x),(x-1,y-1;-y+1,-x+1)}=q^2$.
\end{itemize}
Note that these are also the values one would find by applying 
Serganova's algorithm~\cite{Serganova2} without the specialization
$q=-1$. 
Using these values for $a_{\la,\mu}(q)$, one can calculate explicitly the
matrix elements of the inverse of $A_q$. In particular, we have
determined the values of $(A_q^{-1})_{0,\mu}$, and found
\[
(A_q^{-1})_{0,(-x,-y;y,x)} = q^{x+y},\qquad y \geq x 
 \geq 0 \qquad (x,y\in\mathbb{Z}),
\]
and 0 elsewhere. This coincides with the known values of $K_{0,\mu}(q)$.
We have constructed only one row of the inverse matrix (namely where we
can compare the answer), 
but this on its own is already a rather strong argument
in favour of the conjecture since determining this row involves
the knowledge of ``all'' elements of $A_q$.

Just as conjecture~\ref{conj1} has some interesting consequences, also
the present conjecture has nice implications, in particular~:
\[
\sum_\la a_{\la,\mu}(q) = (1-q)^{\#\mu}
\]
for every $\mu\in P^+$. And therefore also~:
\[
\sum_\la K_{\la,\mu}(q) = {1 \over (1-q)^{\#\mu} }.
\]

\section*{Acknowledgements}
The authors would like to thank V.\ Serganova for her interest and
for useful comments.
J.\ Van der Jeugt wishes to thank the Department 
of Pure Mathematics of the University of Adelaide and in particular
R.B.\ Zhang for their kind hospitality during 
his stay in Adelaide (May--June 1998)
when most of this work was done.

\end{document}